\input amstex
\documentstyle{amsppt}
\refstyle{C}
\TagsOnRight
\define\e{@,@,@,@,@,}
\define\ng{\negthickspace}
\define\ee{\varepsilon_{\chi}}
\define\epo{\varepsilon_{\chi^0}}
\define\epi{\varepsilon_{\chi^i}}
\define\ig{\sigma}
\define\hox{\widehat{H}_{\chi}^{0}}
\define\hoxi{\widehat{H}_{\chi^{i}}^{0}}
\define\mx{M^{\e\chi}}
\define\cok{\roman{Coker}\,}
\define\ra{\rightarrow}
\define\krn{\roman{Ker}\,}
\define\be{@!@!@!@!@!}
\font\tencyr=wncyr10
\font\sevencyr=wncyr7
\font\fivecyr=wncyr5
\newfam\cyrfam
\def\cyr{\fam\cyrfam\tencyr}
\textfont\cyrfam=\tencyr
\scriptfont\cyrfam=\sevencyr
\scriptscriptfont\cyrfam=\fivecyr
\loadbold
\date{April 16th, 2002.}\enddate
\topmatter
\title Finite modules over non-semisimple group rings\endtitle
\author Cristian D. Gonzalez-Aviles\endauthor
\affil{Facultad de Ciencias, Universidad de
Chile}\endaffil
\address{Casilla 653, Santiago, Chile.}
\endaddress
\email{cgonzale\@uchile.cl}\endemail

\bigskip

\thanks{Supported by Fondecyt grant 1000814.}\endthanks

\abstract{Let $G$ be an abelian group of order $n$ and let $R$ be a commutative ring which admits a homomorphism ${\Bbb Z}[\zeta_{n}]\ra R$, where $\zeta_{n}$ is a
(complex) primitive $n$-th root of unity. Given a finite $R[G\e]$-module $M$, we derive a formula relating the order of $M$ to the product of the orders of the various isotypic components $M^{\e\chi}$ of $M$, where $\chi$ ranges over the group of $R$-valued characters of $G$. We then give conditions under which the order of $M$ is exactly equal to the product of the orders of the $M^{\chi}$. To derive these conditions, we build on work of E.Aljadeff and obtain, as a by-product of our considerations, a new criterion for cohomological triviality which improves the well-known criterion of T.Nakayama. We also give applications to abelian varieties and to class groups of abelian fields, obtaining in particular some new class number formulas. Our results also have applications to ``non-semisimple" Iwasawa theory, but we do not develop these here. In general, the results of this paper can be used to strengthen a variety of known results involving finite $R[G\e]$-modules whose hypotheses include (an equivalent form of) the following assumption: ``the order of $G$ is invertible in $R$".}

\endabstract

\endtopmatter

\heading 0. Introduction\endheading
Let $A$ be an abelian variety defined over a global field $F$ and let $K/F$ be a quadratic extension
with Galois group $G$. Write $A^{\roman{t}}$ for the abelian variety dual to $A$. For each of the two characters $\chi$ of $G$, let
$A^{\chi}$ (resp. $(A^{\roman{t}})^{\chi}$) be the $\chi$-twist of $A$ (resp. $A^{\roman{t}}$). In [6] the following result was established.
\proclaim\nofrills{Theorem 0.1\,}{\rm([6], Corollary 4.6)} With the above notations, assume that the following conditions hold.
\roster
\item"(i)" The Tate-Shafarevich group ${\cyr W}(A_{K})$ of $A_{K}$ is finite.
\item"(ii)" $\widehat{H}^{i}(G,A^{\chi}(K))=\widehat{H}^{i}(G, (A^{\roman{t}})^{\chi}(K))=0$ for all integers $i$ and all characters $\chi$ of $G$.
\item"(iii)" Both $A(F_{v})$ and $A^{\roman{t}}(F_{v})$ are connected for all real primes $v$ of $F$.
\endroster
Then
$$[{\cyr W}(A_{K})]=[{\cyr W}(A_{F})]\e
[{\cyr W}(A^{\chi}_{\be\be F})]\cdot\prod_{v\in T}
[H^{1}(G_{w},A(K_{w}))],$$
where $T$ is the set of primes of $F$ which ramify in $K/F$ or where $A_{F}$ has bad reduction, and, for each $v\in T$, $w$ is a fixed prime of $K$ lying above $v$ and $G_{w}=\roman{Gal}(K_{w}/F_{v})$.
\endproclaim

In our attempts to generalize the above theorem to extensions $K/F$ of degree greater than 2, we were
led to the following general problem. Given a finite abelian group $G$ and finite $R[G\e]$-module $M$ (where $R$ is a commutative ring which contains the values of all characters of $G$, e.g. $R=\Bbb Z[\zeta_{n}]$, where $\zeta_{n}$ is a complex $n$-th root of unity), find a formula for the order of $M$ in terms of the orders of the various isotypic components $M^{\chi}$ of $M$, where $\chi$ runs over the group of characters $\widehat{G}$
of $G$ and $M^{\chi}=\{m\in M\: \ig\e m=\chi(\ig)m\,\,\text{for all }\ig\in G\}$. If $n$ denotes the order of $G$ and one considers $R_{*}=\Bbb Z[1/n]\otimes_{\Bbb Z} R$, then it is quite easy to find a formula of the desired type for the order of $M_{*}=\Bbb Z[1/n]\otimes_{\Bbb Z}M$ (which we regard as an $R_{*}[G\e]$-module in the natural way), because there is an isomorphism
$$M_{*}\simeq\bigoplus_{\chi} M_{*}^{\e\chi}=\bigoplus_{\chi}\e\varepsilon_{\chi}^{\prime}M_{*},$$
where $\varepsilon_{\chi}^{\prime}=(1/n)\otimes
\sum_{\sigma\in G}\overline{\chi}(\ig)\e\ig$ is the idempotent of the group ring $R_{*}[G\e]$ corresponding to $\chi\in\widehat{G}$ (here $\overline{\chi}$ denotes the inverse of $\chi$). However, if (for example) $nM=0$, then $M_{*}=0$ and no information is gained on the order of $M$. A different approach
involves the ``quasi-idempotents"
$$\varepsilon_{\chi}=\sum_{\sigma\in G}\overline
{\chi}(\ig)\e\ig\,\,\,\in R[G]$$
(note that these elements satisfy $\varepsilon_{\chi}^{2}=n\e\varepsilon_{\chi}$).
The first observation that we make is that $\varepsilon_{\chi}M$ is no longer equal to $M^{\e\chi}$, so
it is reasonable to expect that the modules
$M^{\e\chi}/\varepsilon_{\chi}M$ will play a role in this work. For $\chi=\chi^{0}$ (the trivial character of $G$), $M^{\e\chi}/\varepsilon_{\chi}M$ is the familiar Tate cohomology module $M^{\e G}/\e N_{G}M=
\widehat{H}^{0}(G, M)$, where $N_{G}=
\sum_{\sigma\in G}\ig$ is the norm element of $R[G\e]$.
In general, $\widehat{H}^{0}_{\chi}(G, M)\overset{\text{def}}\to{=}M^{\e\chi}/\varepsilon_{\chi}M=\widehat{H}^{0}(G, M_{\e\overline{\chi}})$, where
$M_{\e\overline{\chi}}$ is the $\overline{\chi}$-twist
of $M$ (see Section 1 for definitions). Then the following holds.
 
\proclaim{Theorem 0.2} If $G$ is cyclic and $M$ is a finite $R[G\e]$-module, then
$$[M\e]\cdot\prod_{\chi\in\widehat{G}}\left[
\widehat{H}^0_{\chi}(G,M)/S_{\chi}(M)\right]=\prod_{
\chi\in\widehat{G}}\e[M^{\e\chi}\e],$$
where $S_{\chi}(M)$ (for $\chi\in\widehat{G}$) is the
submodule of $\widehat{H}^0_{\chi}(G,M)$ defined in Section 2, formula (6). Equivalently,
$$[M\e]=\prod_{\chi\in\widehat{G}}\left[\e\e\varepsilon_{\chi} M\e\right]\cdot\prod_{\chi\in\widehat{G}}
\left[\e S_{\chi}(M)\e\right].$$
\endproclaim
(Similar formulas exist for any finite abelian group $G$. See Theorem 3.1 below.)

If $M$ is $n$-divisible, i.e., has no $n$-torsion, then the modules $\widehat{H}^0_{\chi}(G,M)$
(and therefore also the modules $S_{\chi}(M)$) vanish. In this case the formulas of the theorem read
$$[M\e]=\prod_{\chi\in\widehat{G}}
\e[M^{\e\chi}\e]=\prod_{\chi\in\widehat{G}}\e
\left[\e\e\varepsilon_{\chi}M\e\right],\tag 1$$
which, by the way, we could also have obtained using the arguments given before the statement of the theorem. The interest of Theorem 0.2 is that there may exist other instances (besides that in which $M$ is $n$-divisible) where the modules $\widehat{H}^0_{\chi}(G,M)=\widehat{H}^0(G,M_{\e\overline{\chi}})$ vanish, and therefore (1) holds. As regards the vanishing of
$\widehat{H}^0(G,M)$ for an arbitrary $G$-module $M\,$
\footnote{Here ``$G$-module" means $\Bbb Z[G\e]$-module.}, in Section 4 we use (a
slightly generalized version of) an impressive theorem of E.Aljadeff [1] and establish the following result.

\proclaim{Theorem 0.3} Let $G$ be any finite group and let $M$ be any $G$-module. Assume that $\widehat{H}^0(H,M)=0$ for all subgroups $H$ of $G$ of {\rm{prime}} order. Then $\widehat{H}^0(H,M)=0$ for all subgroups $H$ of $G$.
\endproclaim
(For an account of the rather unusual and most fortunate way in which I became aware of the
existence of Aljadeff's result, see below.)

Still in Section 4, we use the duality theorem for finite groups to deduce from the previous theorem
the following striking criterion for cohomological triviality.

\proclaim{Theorem 0.4}
Let $G$ be a finite group and let $M$ be a $G$-module. Assume that $\widehat{H}^{-1}(H,M)=\widehat{H}^{0}(H,M)=0$
for all subgroups $H$ of $G$ of prime order. Then $M$ is cohomologically trivial.
\endproclaim
(We remind the reader that a $G$-module $M$ is said to be {\it cohomologically trivial} if $\widehat{H}^{i}(H,M)=0$ for all integers $i$ and all subgroups $H$ of $G$.)

The above criterion is a significant improvement of the well-known criterion established by
T.Nakayama in the mid 1950's (see [2,\S 9]). For further comments on the above result, see the remarks following the proof of Theorem 4.6 below.

Using the above criterion, we obtain results of the following type.

\proclaim{Theorem 0.5} Let $p$ be a prime number
and let $G$ be a cyclic group of order $p^n$,
where $n\geq 1$. Let $H$ be the unique subgroup
of $G$ of order $p$ and let $M$ be a finite $R[G\e]$-module. Assume that $\widehat{H}_{\chi}^{0}(H,M)=0$ for all characters $\chi$ of $H$. Then
$\widehat{H}_{\chi}^{0}(G,M)=0$ for all characters $\chi$ of $\, G$ and
$$[M\e]=\prod_{\chi\in\widehat{G}}[M^{\e\chi}\e]=
\prod_{\chi\in\widehat{G}}[\e\varepsilon_{\chi} M\e].$$
\endproclaim

\proclaim{Corollary 0.6} Let $G$ be a cyclic group of order $2^n$, where $n\geq 1$, and let $\tau$ be
the unique element of $G$ of order 2. Let $M$ be a
finite $R[G\e]$-module such that $M^{+}\overset{\text{def}}\to{=}\{m\in M\:
\tau\e m=m\}=(1+\tau\e)M$. Then
$$[M\e]=\prod_{\chi\in\widehat{G}}[M^{\e\chi}\e]
=\prod_{\chi\in\widehat{G}}[\e\varepsilon_{\chi}\be M\e].$$
\endproclaim

In Section 5 we apply the results of the preceding sections to Tate-Shafarevich groups of abelian varieties and obtain, in particular, the following generalization of Theorem 0.1.

\proclaim{Theorem 0.7} Let $K/F$ be a cyclic Galois
extension of global fields, with Galois group $G$ of order $n$, and let $A$ be an abelian variety defined over
$F$ with complex multiplication by the ring of integers of $\,{\Bbb Q}(\zeta_{n})$. Assume that
the following conditions hold.
\roster
\item"(i)" The Tate-Shafarevich group ${\cyr W}(A_{K})$ of $A_{K}$ is finite.
\item"(ii)" $\widehat{H}_{\chi}^{i}(G,A(K))=
\widehat{H}_{\chi}^{i}(G,A^{\roman{t}}(K))=
0$ for all $i$ and all $\chi\in\widehat{G}$.
\item"(iii)" $A(F_{v})$ is connected for all real primes $v$ of $F$.
\endroster
Then
$$[{\cyr W}(A_{K})]=\prod_{\chi\in\widehat{G}}
[{\cyr W}(A^{\chi}_{\be\be F})]\cdot\prod_{\chi\in\widehat{G}}[S_{\chi}(
{\cyr W}(A_{K}))],$$
where, for each $\chi\in
\widehat{G}$, $S_{\chi}({\cyr W}(A_{K}))$ is the submodule of $\widehat{H}_{\chi}^{0}(G,{\cyr W}(A_{K}))$ defined in Section 2, formula (6).
Further, for each $\chi\in\widehat{G}$, the order of $S_{\chi}({\cyr W}(A_{K}))$ divides 
$$\prod_{v\in T_{\e\overline{\chi}}}[\widehat{H}_{\chi_{_{w}}}^{0}\be\be(G_{w},A(K_{w}))],$$
where $T_{\e\overline{\chi}}$ and $\chi_{_{w}}$ are as in the statement of Corollary 5.6 below.
\endproclaim

In Section 6 we apply the results of Sections 2-4 to study class groups of abelian fields and obtain the following results.

\proclaim{Theorem 0.8} Let $K/F$ be a finite abelian
extension of exponent 2. Then there exists an integer $t$ such that
$$2^{t}\cdot h_{K}/h_{F}=\prod_{\chi\neq\chi^{0}}
\e[\varepsilon_{\chi}C_{K}],$$
where $C_{K}$ is the ideal class group of $K$,
$h_{K}$ (resp. $h_{F}$) denotes the order of $C_{K}$ (resp. $C_{F}$) and the product extends over all non-trivial characters of $G$.
\endproclaim
(Regarding the above result, there are indications that the integers $[\varepsilon_{\chi}C_{K}]$ appearing on the right-hand side of the above formula are related to the class numbers of the various subextensions of $K/F$. See the remarks following the statement of Theorem 6.2 below, where we also comment on the probable value of $t$.)

If $K/F$ does not have exponent 2, then it is necessary to extend scalars. Let $\overline{C}_{K}=\Bbb Z[\zeta_{n}]
\otimes_{\Bbb Z}C_{K}$, where $n$
is the degree of $K/F$ and $\zeta_{n}$ is a
fixed complex $n$-th root of unity. Let
$\varphi(n)=[(\Bbb Z/n\Bbb Z)^{\times}]$. Then

\proclaim{Theorem 0.9} Let $K/F$ be a cyclic Galois extension of number fields with Galois group $G$ of order $n$. Assume that $K/F$ contains no subextensions which are everywhere-unramified. Then
$$\left(h_{K}/h_{F}\right)^{\varphi(n)}=\prod_{\Sb \chi\in\widehat{G}\\\chi\neq\chi^{0}\endSb}
\e[S_{\chi}(\overline{C}_{K})]\cdot
\prod_{\Sb \chi\in\widehat{G}\\\chi\neq\chi^{0}\endSb}
\e[\varepsilon_{\chi}\overline{C}_{K}],$$
where the integers $[S_{\chi}(\overline{C}_{K})]$ are divisible only by the primes that divide $n$.
\endproclaim

\proclaim{Theorem 0.10} Let $K/F$ be a cyclic Galois extension with Galois group $G$ of order $2^{n}$, where $n\geq 1$. Assume that $K/F$ ramifies at some prime. Then
$$(h_{K}/h_{F})^{2^{n-1}}=\prod_
{\chi\neq\chi^{0}}[\e\varepsilon_{\chi}\be
\overline{C}_{K}].$$
\endproclaim

\proclaim\nofrills{Corollary 0.11\,}{\rm(of Theorem 0.9)}. Let $p$ be an odd prime and let $K=\Bbb Q(\zeta_{p})^{+}$ be the maximal real subfield of $\Bbb Q(\zeta_{p})$. Write $h^{+}$ for the class number of $K$ and $G$ for the Galois group of $K/\Bbb Q$. Then
$$(h^{+})^{\varphi\left(\frac{p-1}{2}\right)}=r\cdot\ng
\prod_{\chi\in\widehat{G}}\e[\varepsilon_{\chi}\overline{C}_{K}],$$
where $r$ is an integer which is divisible only by the primes that divide $(p-1)/2$.
In particular, $p$ divides $h^{+}$ if and only if $p$ divides $[\varepsilon_{\chi}\overline{C}_{K}]$ for
some character $\chi$ of $G$.
\endproclaim

The above results (0.8-0.11) cannot be considered satisfactory, because they give no information on the integers $[\varepsilon_{\chi}\overline{C}_{K}]$. 
As mentioned above, these integers seem to be related to the class numbers of the various subextensions
of $K/F$, a conjecture which we hope to confirm in a future publication.
As regards the last assertion of Corollary 0.11, it is of course an allusion to Vandiver's conjecture, which asserts that  $h^{+}$ is never divisible by $p$. Regarding this conjecture, the results of this paper seem to indicate that the following statement is true: Vandiver's conjecture holds for $p$ if and only if $p$ does not divide $h_{L}$ for every subextension $L/\Bbb Q$ of $\Bbb Q(\zeta_{p})^{+}/\e\Bbb Q$ of {\it prime} degree. See Section 6 for additional comments.

We also note that the results of this paper have applications to ``non-semisimple" Iwasawa theory, which we hope to develop in a future publication. 
\smallpagebreak

Finally, I would like to share with the reader the following anecdote connected with the writing of this paper. In August of 2001 I attended the {\it XIV Coloquio Latinoamericano de Algebra}, which took place in La Falda, Argentina. There I had the opportunity to meet many distinguished mathematicians, including V.Kac, B.Kostant and D.Sullivan. At one point during the conference I noticed D. Sullivan and E. Aljadeff having a conversation in the main hall. After observing them for a few seconds, it became clear to me that Sullivan wanted to part ways in order to attend to some other bussiness, but did not want to leave Aljadeff alone. So as I was passing by them Sullivan waved to me to stop and said to Aljadeff: ``Eli, tell Chris what you just told me", and then left. Aljadeff then proceeded to explain to me his result on the surjectivity of the norm map (see Theorem 4.1 below), which was exactly the type of result I needed to complete this paper. The conclusion is clear. D.Sullivan furthers the advancement of Mathematics even without trying.

\bigskip

\heading Acknowledgements \endheading
I am very grateful to E. Aljadeff for explaining to me his results on the norm map and for providing me with a copy of his paper [1]. {\it Je remercie \'egalement l'Officine pour la Cooperation R\'egionale de l'Ambassade de France au Chili, et tr\`es particuli\`erement Mme. Annie Pauget, pour le financement de mon voyage \`a Sierra de La Ventana, Argentine,  o\'u une version pr\'eliminaire de ce travail a et\'e present\'e le 30 d'Octobre 2001, dans le cadre du} Coloquio Internacional de Homolog\'{\i}a y Representaciones de Algebras.

\heading 1. Preliminaries\endheading

Let $G$ be a finite group of exponent $e$ and
let $\chi\: G\ra F$ be a character (one-dimensional representation) of $G$ with values in some field $F$ which contains a primitive
$e$-th root of unity $\zeta_e$ (e.g. $F=\Bbb C$).
A commutative ring (with unit) $R$ is said to be {\it large enough} for $G$ if there exists a ring homomorphism $\Bbb Z[\zeta_e]\ra R$. For example, $\Bbb Z[\zeta_e]$ and ${\Bbb F}_{p}[x]/(x^{p^m})$
($e=p^m$ a prime power) are large enough for $G$. Since $\chi\: G\ra F$ factors through $Z[\zeta_e]$, we may compose $\chi$ with the given homomorphism $Z[\zeta_e]\ra R$ to obtain a multiplicative map $G\ra R$, which will also be denoted by $\chi$. Thus $\chi\: G\ra R$ is a ``character of $G$ with values in $R$". Clearly the values of $\chi$ lie in $R^{\times}$ (the group of units of $R$), so $\chi$ is an element of $\widehat{G}=\roman{Hom}(G,R^{\times})$.

Now let $M$ be an $R[G\e]$-module, where $R$ is large enough for $G$, and consider the augmentation
homomorphism
$$\alpha_{\chi}\: R[G\e]\ra R,\quad
\sum_{\ig\in G}r_{\ig}\e\ig\mapsto\sum_{\ig\in G}r_{\ig}\e\chi(\ig).$$
For each $i\geq 0$, we define the $R$-modules
$$H_{\chi}^i(G,M)=\roman{Ext}_{R[G\e]}^i(R,M),$$
where $R$ is being regarded as an $R[G\e]$-module via the map $\alpha_{\chi}$. For $i=0$ we have
$$H_{\chi}^0(G,M)=\roman{Hom}_{R[G\e]}(R,M)=M^{\e\chi},$$
where $M^{\e\chi}=\{m\in M\: \sigma\e m=\chi(\sigma)\e m\quad\text{for all }\ig\in G\}$.

Now consider the element
$$\ee=\sum_{\ig\in G}\overline{\chi}(\ig)\e\ig\quad
\in R[G\e],$$
where $\overline{\chi}$ denotes the inverse of $\chi$ (i.e., $\overline{\chi}(\sigma)=\chi(\sigma)^{-1}$ for all $\sigma\in G$). It is not difficult to check that $\ee M\subset\mx$. Define
$$\hox(G,M)=\mx/\ee M.$$
If $M$ is a ${\Bbb Z}[G\e]$-module and $\chi$ is the trivial character of $G$, then the $\Bbb Z$-modules 
$H_{\chi}^i(G,M)$ and $\hox(G,M)$ defined above are the well-known groups $H^0(G,M)=M^{\e G}$ and $\widehat{H}^0(G,M)=M^{\e G}/N_{\be G}M$ from Group Cohomology. In general, the $R$-modules
$H_{\chi}^i(G,M)$ and $\hox(G,M)$ are ``usual"
cohomology modules of a twisted form of $M$. Indeed,
if $M_{\e\overline{\chi}}$ denotes the $R$-module $M$
endowed with the new $G$-action
$$\sigma\cdot m=\overline{\chi}(\sigma)\e\sigma\e m
\qquad(\sigma\in G,\; m\in M),$$
then
$$\widehat{H}^i_{\chi}(G,M)=\widehat{H}^i(G,
M_{\e\overline{\chi}})$$
for all $i\geq 0$ (this follows in a standard way
from the fact that the functor $M\mapsto M^{\e\chi}$ is the composite of the functors $M\mapsto M^{\chi^0}$ and $M\mapsto M_{\e\overline{\chi}}$, the second of which is {\it exact} and right adjoint to an exact functor, namely $M\mapsto M_{\chi}$.)

\heading 2. Cyclic groups\endheading
We assume now that $G$ is a finite cyclic group
and write $n$ for its order. Let $M$ be a finite $R[G\e]$-module, where $R$ is large enough for $G$. Then there exists a ring homomorphism $\Bbb Z[\zeta_{n}]\ra R$, where $\zeta_{n}$ is a primitive (say complex) $n$-th root of unity. We will continue to write $\zeta_{n}$ for the image of 
$\zeta_{n}$ in $R$ under the above homomorphism (this should cause no confussion). We now choose and fix a generator $\tau$ of $G$ and define a character
$\chi\: G\ra R^{\times}$ by $\chi(\tau)=\zeta_{n}$.
Then any other character of $G$ is of the form $\chi^i$, where $0\leq i\leq n-1$ (by convention $\chi^0$ is the trivial character of $G$, i.e.,
$\chi^0(\sigma)=1$ for all $\sigma\in G$).

Now recall the elements
$\epi=\sum_{\ig\in G}\overline{\chi}^{\e i}(\ig)\e\ig\,\in R[G\e]$, where $0\leq i\leq n-1$. We have
$$\align\epi =\sum_{j=0}^{n-1}\e\overline{\zeta}_{n}^{\e\e ij}
{\tau}^{j}&=\prod_{j=1}^{n-1}(\e\overline{\zeta}_{n}^{\e\e i}\e\tau-\zeta_{n}^{\e j})\\
&=\overline{\zeta}_{n}^{\e\e (n-1)i}\prod_{j=1}^{n-1}(\tau-\zeta_{n}^{\e i+j})=\zeta_{n}^{\e i}\prod_{\Sb j=0\\j\neq i\endSb}^{n-1}(\tau-\zeta_{n}^{\e j}).\endalign$$
In particular the norm element
$\epo=\sum_{\ig\in G}\ig\in R[G\e]$ factors as
$$\epo=\prod_{j=1}^{n-1}(\tau-\zeta_{n}^{j}).$$
Our objective now is to derive a formula
relating the order of $M$ to the orders of the various isotypic components $M^{\chi^i}$ of $M$. We will use the following notation. If $a$ is any element of $R[G\e]$, $\krn a$ will denote the kernel  of multiplication-by-$a$ on $M$ (note that this is
simply a nonstandard notation for the $a$-torsion submodule of $M$). We note that $M^{\e\chi^i}=\{m\in M\: \tau\e m=\zeta_{\e n}^{\e i}\e m\}=\krn(\tau-\zeta_{\e n}^{\e i})$. The order of a finite module $M$ will be denoted by $[M]$. Now, in order to make our general arguments
more transparent, we will begin by examining the simplest case, that in which $n=2$. Consider the following exact sequence (which is available for any $n$)
$$0\ra\krn\epo\ra M\overset{\varphi_{0}}
\to{\longrightarrow}M^{\e\chi^0}\ra\widehat{H}^0_{\chi^0}(
G,M)\ra 0\tag 2$$
where $\varphi_{0}$ is the multiplication-by-$\e\epo$ map. When $n=2$, $\krn\epo=\krn(1+\tau)=\krn(\tau-(-1))=M^{\e\chi}$ (see above), so we immediately get from (2) the identity
$$[M\e]\e\left[\widehat{H}^0_{\chi^0}(G,M)\right
]=[M^{\e\chi^0}\e]\e [M^{\e\chi}],$$
which is the desired result for $n=2$.

When $n=3$ the situation is more complicated, because $\epo=(\tau-\zeta_{3})(\tau-\zeta_{3}^{2})$ is the product of {\it two} linear factors, and therefore $\krn\epo$ cannot equal $M^{\chi^i}$ for any $i$. However, we can relate $\krn\epo$ to
modules of the form $M^{\chi^i}$ by means of the exact sequence
$$0\ra M^{\e\chi^2}\ng\be=\be\krn(\tau-\zeta_{3}^2)\ra\krn\epo
\overset{\varphi_1}\to{\longrightarrow}
M^{\e\chi}\ra Q_1\ra 0,$$
where $\varphi_1$ is the multiplication-by-
$\ng(\tau-\zeta_{3}^{2})$ map and
$$Q_1=\cok\varphi_1= M^{\e\chi}/(\tau-\zeta_{3}^{2})\e\krn\epo.$$
Now since
$$(\tau-\zeta_{3}^{2})\e\krn\epo=M^{\e\chi}\cap(\tau-\zeta_{3}^{2})M,$$
we have
$$Q_1\simeq\widehat{H}_{\chi}^0(G,M)/S_1(M),$$
where
$$S_1(M)=M^{\e\chi}\cap(\tau-\zeta_{3}^{2})M/
{\varepsilon}_{\chi}M.$$
Thus, writing $S_0(M)=\{0\}$,  we have
$$[M\e]\left[\widehat{H}^0_{\chi^0}(G,M)/S_0(M)
\right]\left[\widehat{H}^0_{\chi}(G,M)/S_1(M)
\right]=
[M^{\e\chi^0}\e][M^{\e\chi}][M^{\e\chi^2}
\e],$$
which is the desired result for $n=3$.

We now present the general argument.

For each $i\in\{0, 1,\dots , n-2\}$, there is an exact sequence
$$0\ra\krn\tsize{\prod\limits_{j=i+1}^{n-1}}(\tau-\zeta_{n}^{j})\ra\krn\tsize{\prod\limits_{j=i}^{n-1}}(\tau-\zeta_{n}^{j})\overset{\varphi_i}
\to{\longrightarrow}M^{\chi^i}\ra Q_i\ra 0\tag 3$$
in which $\varphi_i$ is the multiplication-by-
$\ng\prod_{j=i+1}^{n-1}(\tau-\zeta_{n}^{j})$ map and
$$Q_i=\cok\varphi_i = M^{\e\chi^i}\Big/\tsize{\prod\limits_{j=i+1}^{n-1}}(\tau-\zeta_{n}^{j})\e\e\roman{Ker}\e\e
\tsize{\prod\limits_{j=i}^{n-1}}(\tau-\zeta_{n}^{j}).$$
Note that (3) with $i=0$ is precisely (2), because $\prod_{j=0}^{n-1}(\tau-\zeta_{n}^{j})=\tau^n-1=0$
and therefore $\krn\prod_{j=0}^{n-1}(\tau-\zeta_{n}^{j})=M$. Also note that, since
$$\prod_{j=i+1}^{n-1}(\tau-\zeta_{n}^{j})\e\e\roman{Ker}\e\e\prod_{j=i}^{n-1}(\tau-\zeta_{n}^{j})=M^{\e\chi^i}\cap\prod_{j=
i+1}^{n-1}(\tau-\zeta_{n}^{j})\, M,$$
we have
$$Q_i\simeq\hoxi(G,M)/S_i(M),\tag 4$$
where $S_i(M)$ is the submodule of $\hoxi(G,M)$ defined by
$$S_i(M)=\Big(M^{\e\chi^i}\cap\prod_{j=i+1}^{n-1}(\tau-\zeta_{n}^{j})\, M\Big)\Big/\epi M.\tag 5$$
Now, by (3),
$$\frac{[\krn\prod_{j=i}^{n-1}(\tau-\zeta_{n}^{j})]}{
[\krn\prod_{j=i+1}^{n-1}(\tau-\zeta_{n}^{j})]}=\frac{
[M^{\e\chi^i}\e]}{[Q_i]}$$
for $i=0, 1,..., n-2$. Multiplying these equalities together and noting that the product of the left-hand side terms {\it telescopes}, we obtain
$$\frac{[M\e]}{[M^{\e\chi^{n-1}}\e]}=\frac{
\prod_{\e i=0}^{\e n-2}\e [M^{\e\chi^i}\e]}{\prod_{\e i=0}^{\e n-2}\e[Q_i]}.$$
Thus, using (4), we conclude
\proclaim{Theorem 2.1} We have
$$[M\e]\cdot\prod\limits_{i=0}^{n-2}\left[
\hoxi(G,M)/S_i(M)\right]=\prod\limits_{i=0}^{n-1}\e\big[M^{\e\chi^i}\e\big],$$
where $S_i(M)$ ($i=0,1,\dots,n-2$) is given by (5).
\endproclaim
We will now restate the above theorem in a
form which is more suitable for generalization. Write
$\prod_{j=i+1}^{n-1}(\tau-\zeta_{n}^{j})\, M=M$ if
$i=n-1$. Then, for $\psi\in\widehat{G}$, set
$$S_{\psi}(M)=S_i(M)=\Big(M^{\e\chi^i}\cap\prod_{j=i+1}^{n-1}(\tau-\zeta_{n}^{j})\, M\Big)\Big/\epi M\tag 6$$
if $\psi=\chi^{i}$ with $0\leq i\leq n-1$. Then Theorem 2.1 can be restated as follows.
\proclaim{Theorem 2.2} If $G$ is a cyclic group and $M$ is a finite $R[G\e]$-module, then
$$[M\e]\cdot\prod_{\psi\in\widehat{G}}\left[
\widehat{H}^0_{\psi}(G,M)/S_{\psi}(M)\right]=\prod_{
\psi\in\widehat{G}}\e[M^{\e\psi}\e],$$
where $S_{\psi}(M)$ ($\psi\in\widehat{G}$) is given by (6). Equivalently,
$$[M\e]=\prod_{\psi\in\widehat{G}}\left[\e\e\varepsilon_{\psi} M\e\right]\cdot\prod_{\psi\in\widehat{G}}
\left[\e S_{\psi}(M)\e\right].$$
\endproclaim


\heading 3. Abelian groups\endheading
In this section we generalize Theorem 2.2 to arbitrary (finite) abelian groups. We will consider first abelian groups which are the direct product of two cyclic groups.

Let $K_1$ and  $K_2$ be (finite) cyclic groups and let
$G=K_1\times K_2$ be the direct product of $K_1$ and $K_2$. Let $M$ be a finite $R[G\e]$-module, where $R$ is large enough for $G$ (for example $R=\Bbb Z[\zeta_{e}]$, where $e$ is the exponent of $G$ and $\zeta_{e}$ is a fixed complex $e$-th root of unity).
Write $M_{K_i}$ for the $R[G\e]$-module $M$ regarded as an $R[K_i]$-module ($i=1,2$). Note that since $\widehat{G}=\widehat{K}_{1}\times\widehat{K}_{2}$, each character
$\chi\in\widehat{G}$ can be written as $\chi=\chi_{_1}\ng\cdot\chi_{_2}$, where $\chi_{_1}\in\widehat{K}_{1}$ and
$\chi_{_2}\in\widehat{K}_{2}$. We have
$$M^{\e\chi}=(M_{K_1}^{\e\chi_1})_{K_2}^{\e
\chi_2}$$
(observe that $M_{K_1}^{\e\chi_1}$ is naturally an
$R[G\e]$-module, so we can consider the restricted module $(M_{K_1}^{\e\chi_1})_{K_2}$). To ease notation, we will write the above equality simply as ``$M^{\e\chi}=(M^{\e\chi_{_1}})^{\e\chi_
{_2}}$". (The reader should bear in mind that in an expression of the form ``$M^{\e\chi_{_1}}$" (resp. ``$N^{\e\chi_{_2}}$"), $M$ (resp. $N$) is being regarded as a $K_1$-module (resp. $K_2$-module).) We have
$$\prod_{\chi\in\widehat{G}}\e[M^{\e\chi}\e]=
\prod_{\chi_{_1}\ng\in\widehat{K}_{1}}\prod_{\,\chi_{_2}\be\be\in\widehat{K}_{2}}\e[(M^{\e\chi_{_1}})^{\e\chi_{_2}}\e].$$
Now by Theorem 2.2,
$$\prod_{\chi_{_2}\be\be\in\widehat{K}_{2}}\e[(M^{\e\chi_{_1}})^{\e\chi_{_2}}\e]=[M^{\e\chi_{_1}}\e]
\cdot\prod_{\chi_{_2}\be\be\in\widehat{K}_{2}}\left[\widehat{H}^0_{\chi_{_2}}\be(K_2,M^{\e\chi_{_1}}\be)/
S_{\chi_{_2}}\be(M^{\e\chi_{_1}}\be)
\right],$$
where $S_{\chi_{_2}}\be\be(M^{\e\chi_{_1}}\be)$ is given by (6) for $\psi=\chi_{_2}$ and $M=M^{\e\chi_{_1}}$. Applying Theorem 2.2 once again, we obtain the formula
$$[M\e]\be\cdot\prod_{\chi_{_1}\ng\in\widehat{K}_{1}}\left[\widehat{H}^0_{\chi_{_1}}\be(K_1,M)
/S_{\chi_{_1}}\be(M)\right]\ng
\prod_{\Sb \chi_{_2}\be\be\in\widehat{K}_{2}\\
\chi_{_1}\ng\in\widehat{K}_{1}\endSb}
\left[\widehat{H}^0_{
\chi_{_2}}\be(K_2,M^{\e\chi_{_1}}\be)/
S_{\chi_{_2}}\be\be(M^{\e\chi_{_1}}\be)
\right]
=\prod_{\chi\in\widehat{G}}\e[M^{\e\chi}\e].$$
In general, the following holds.
\proclaim{Theorem 3.1} Let $G$ be a finite abelian group and let $M$ be a finite $R[G\e]$-module. Suppose $G=K_1\times\dots\times K_r$ ($r\geq 1$) is
a decomposition of $G$ as a direct product of cyclic groups. For $0\leq i\leq r-1$, set $G_{i}=K_1\times\dots\times K_{i}
$, where $G_{0}$ is defined to be $\{0\}$. Then
$$[M\e]\cdot\prod_{i=1}^{r}\prod_
{\Sb \chi_{_i}\be\be\in\widehat{K}_{i}\\
\psi\in\widehat{G}_{i-1}\endSb}
\ng\ng\left[\widehat{H}^0_{\chi_{_i}}(K_i,M^{\e\psi}\e)/S_{\chi
_{_i}}(M^{\e\psi}\e)\right]=\prod_{
\chi\in\widehat{G}}\e[M^{\e\chi}\e].$$
\endproclaim
{\it Proof}. This can be proved easily by induction,
writing $G=(K_1\times\dots K_{r-1})\times K_r$ and proceeding as in the case $r=2$.\qed

{\it Remark}. We report that we were unable to find a proper generalization of the results of this section to the case of non-abelian groups $G$. This problem seems to require tools from the representation theory of finite groups, with which we are unfamiliar at the present time. We hope to return to this issue in the future.

\heading 4. A criterion for cohomological
triviality.\endheading
In this section we establish a new criterion for cohomological triviality using an impressive theorem of E.Aljadeff. We then apply this criterion to derive sufficient conditions under which the order of a finite $R[G\e]$-module $M$ (where $G$ is abelian and $R$ is large enough for $G$) equals the product of the orders of the various isotypic components $M^{\e\chi}$ of $M$, where $\chi$ ranges over $\widehat{G}$. In this section, ``$G$-module" means $\Bbb Z[G\e]$-module.

We begin by recalling Aljadeff's result.
\proclaim\nofrills{Theorem 4.1\,}{\rm (Aljadeff)}. Let $G$ be a finite group and let $K$ be a $G$-module which is also
a commutative ring with unit. Suppose that $\widehat{H}^0(H,K)=0$ for all subgroups $H$ of $G$ of \,{\rm{prime}}\, order. Then
$\widehat{H}^0(H,K)=0$ for all subgroups $H$ of $G$.
\endproclaim
{\it Proof}. See [1, Theorem 0.1].\qed
\smallpagebreak

{\it Remark}. Aljadeff's proof of the above result uses certain modules over the twisted group ring $K_{\e\roman{t}}\e G$, where $t\: G\ra\roman{Aut}(K)$ is the map
inducing the action of $G$ on $K$. An alternative proof of the theorem, for $G$ {\it abelian}, can be given (with some effort) using the inflation-restriction exact sequence in group cohomology together with [2, \S 9, Lemma 1].
\smallpagebreak

We will now extend Aljadeff's result to arbitrary $G$-modules. Note that
any $G$-module $M$ is in
particular an abelian group, so it is meaningful to consider the group ring $\Bbb Z[M]$, which is a $G$-module (with $G$ acting trivially on $\Bbb Z$) as well as a commutative ring with unit.
\proclaim{Proposition 4.2} Let $H$ be a group of prime order and let $M$ be an $H$-module. Write $N_{\be H}$ for the norm element of $\Bbb Z[H]$. Then
$$\Bbb Z[M]^{H}=\Bbb Z\left[M^{H}\right]+N_{\be H}\e\Bbb Z[M].$$
\endproclaim
{\it Proof}. To prove the nontrivial containment, we need only check that $\Bbb Z[M\setminus M^{H}]^{H}\subset N_{\be H}\e\Bbb Z[M]$. Let $\sum_{m\in
M\setminus M^{H}}\e a_m\e m\,$ be an element of $\Bbb Z[M\setminus M^{H}]^{H}$, with $a_m=0$ for all but finitely many $m\in M
\setminus M^{H}$. Then
$a_m=a_{\sigma m}$ for all $m\in M\setminus M^{H}$ and any $\sigma\in H$. Thus we can write
$$\sum_{m\in M\setminus M^{H}}\e a_m\e m=\sum_{i}\e
a_i\sum_{m\in M_i}m,$$
where the $M_i$ are the various orbits of the action of $H$ on $M\setminus M^{H}$ and $a_i$ is the common value of $a_m$ as $m$ ranges over $M_i$. Now the stabilizer in $H$ of an element $m\in M\setminus M^{H}$ is necessarily trivial, because it is a proper subgroup of a group of prime order. It follows that $\sum_{m\in M_i}m=N_{\be H}\e m_i$ for any $m_i\in M_i$, and the proposition follows.\qed

\proclaim{Corollary 4.3} Let $H$ be a group of prime order and let $M$ be an $H$-module. Suppose that $\widehat{H}^0(H,M)=0$. Then
$\widehat{H}^0(H,\Bbb Z[M])=0$.
\endproclaim

\proclaim{Theorem 4.4} Let $G$ be any finite group and let $M$ be any $G$-module. Suppose that $\widehat{H}^0(H,M)=0$ for all subgroups $H$ of $G$ of prime order. Then $\widehat{H}^0(H,M)=0$ for all subgroups $H$ of $G$.
\endproclaim
{\it Proof}. Let $H$ be any subgroup of $G$ and consider the surjective $H$-module homomorphism $\Bbb Z[M]\ra M$ defined by mapping $\sum a_m\e m\in\Bbb Z[M]$ to the actual sum in $M$. This map admits a section as a map of $H$-modules, namely $m\mapsto 1\be\cdot m$. Therefore
$\widehat{H}^0(H,M)$ is a direct summand of $\widehat{H}^0(H,\Bbb Z[M])$. On the other hand, by hypothesis and Theorem 4.1 together with Corollary 4.3, $\widehat{H}^0(H,\Bbb Z[M])=0$ for every subgroup $H$ of $G$. The theorem is now clear.\qed

Next, we will use Theorem 4.4 and the duality theorem for finite groups to
establish the criterion for cohomological triviality mentioned earlier.

For any abelian group $M$, we will write $M^{*}=\roman{Hom}(M,\Bbb Q/\Bbb Z)$ for the Pontrjagyn dual of $M$.
\proclaim{Lemma 4.5} Let $G$ be a finite group and let $M$ be any $G$-module. Suppose that $\widehat{H}^{-1}(H,M)=0$ for all subgroups $H$ of $G$ of prime order. Then $\widehat{H}^{-1}(H,M)=0$ for all subgroups $H$ of $G$.
\endproclaim
{\it Proof}. By the duality theorem for finite groups [3, p.250], for every subgroup $H$ of $G$ of prime order
$$\widehat{H}^0(H,M^{*})\simeq\widehat{H}^{-1}(H,M)^{*}=0.$$
Consequently, by Theorem 4.4, $\widehat{H}^0(H,M^{*})=0$
for every subgroup $H$ of $G$. Applying the duality theorem once again, we conclude that $\widehat{H}^{-1}(H,M)=0$ for every $H$.\qed

Recall that, if $G$ is a finite group, a $G$-module $M$ is said to be {\it cohomologically trivial} if
$\widehat{H}^{i}(H,M)=0$ for all subgroups $H$ of $G$ and all integers $i$.
\proclaim{Theorem 4.6} Let $G$ be a finite group and let $M$ be a $G$-module. Assume that $\widehat{H}^{-1}(H,M)=\widehat{H}^{0}(H,M)=0$ for all subgroups $H$ of $G$ of prime order. Then $M$ is cohomologically trivial.
\endproclaim
{\it Proof}. Let $p$ be a prime number and let $G_p$ be a $p$-Sylow subgroup of $G$. Then, by Theorem 4.4 and Lemma 4.5, $\widehat{H}^{-1}(G_p,M)=\widehat{H}^{0}(G_p,M)=0$. Now
the well-known criterion of T.Nakayama [2, \S 9, Theorem 9(i)] completes the
proof.\qed
\smallpagebreak

{\it Remarks}. (a) By the periodicity of the cohomology of cyclic groups [2], the hypothesis of the theorem can be replaced by: $\widehat{H}^{i}(H,M)=\widehat{H}^{j}(H,M)=0$ for all subgroups $H$ of $G$ of prime order and any pair of integers $i,\e j$  which differ by an odd integer.

(b) Theorem 4.6 is a significant improvement of the criterion of T.Nakayama used in its proof. For example, if $G$ is a cyclic group of order $p^{n}$ (where $p$ is a prime and $n\geq 1$) and $M$ is a $G$-module, then the cohomological triviality of $M$ can be checked ``at the first layer", i.e., by checking whether $\widehat{H}^{i}(H,M)$ is zero for $i=-1, 0$ and $H$ the unique subgroup of $G$ of order $p$. Nakayama's criterion requires that one check the vanishing of these cohomology groups for the full group (i.e., for $H=G$), a verification that depends on the (possibly very large) value of $n$.

(c) We also note that Theorem 4.6 can be used to weaken the hypotheses of some well-known theorems in group cohomology which
depend on Nakayama's criterion, e.g., the Nakayama-Tate theorem.

\smallpagebreak

Recall that if $H$ is a finite cyclic group and $M$ is an $H$-module the groups
$\widehat{H}^{i}(H,M)$ for $i=-1,0$ are given by
$$\widehat{H}^{-1}(H,M)=\sideset_{N_{\be H}}\and\to{\be\be M}/I_{\be H} M\quad\text{and}
\quad\widehat{H}^{0}(H,M)=M^{H}/N_{\be H} M,$$ where $\sideset_{N_{\be H}}\and
\to{\be\be M}$ is (the standard notation for) the kernel of multiplication by $N_{\be H}$ on $M$ and $I_{\be H}$ is the augmentation ideal of $\Bbb Z[H]$ (see for example [2]). Thus Theorem 4.6 can be restated as
\proclaim{Theorem 4.7} Let $G$ be a finite group and let $M$ be
any $G$-module. Assume that for each subgroup $H$ of $G$ of prime order we have
$$\sideset_{N_{\be H}}\and\to{\be\be M}=I_{\be H}\e M\qquad\text{and}\qquad M^{H}=N_{\be H}\e M.$$
Then $M$ is cohomologically trivial.\qed
\endproclaim

Recall that for any subgroup $H$ of $G$ (where $G$ is a finite group) and any $G$-module $M$, $M_{H}$ denotes the
abelian group $M$ regarded as an $H$-module. Recall also that
if $H$ is cyclic and the groups
$\widehat{H}^{i}(H,M)$ ($i=0, 1$) are {\it finite}, then the Hebrand quotient $h(M)$ of $M_{H}$ is defined by
$$h(M)=[\widehat{H}^{0}(H,M)]/[H^{1}(H,M)].$$
A well-known theorem states that the Herbrand quotient of a finite module is 1.

\proclaim{Corollary 4.8} Let $G$ be a finite group and let $M$ be a $G$-module.
Assume that for each subgroup $H$ of $G$ of prime order the Herbrand quotient of $M_H$ is defined and equal to 1 (for example, $M$ can be a finite $G$-module). Assume further that $\widehat{H}^{0}(H,M)=0$ for all subgroups $H$ as above. Then $M$ is cohomologically trivial.
\endproclaim
\smallpagebreak

{\it Remark}. The above corollary will be applied in the next section to obtain an interesting result concerning abelian varieties defined over local fields. See Theorem 5.3 below.

\smallpagebreak

Next, we will combine Theorems 3.1 and 4.4 to derive sufficient conditions under which the order of a finite $R[G\e]$-module $M$,
where $G$ is a finite abelian group and $R$ is large enough for $G$, equals the product of the orders of the various isotypic components $M^{\e\chi}$ of $M$ as $\chi$ ranges over $\widehat{G}$.

We begin with
\proclaim{Theorem 4.10} Let $p$ be a prime number
and let $G$ be a cyclic group of order $p^n$,
where $n\geq 1$. Let $H$ be the unique subgroup
of $G$ of order $p$ and let $M$ be a finite $R[G\e]$-module. Assume that $\widehat{H}_{\chi}^{0}(H,M)=0$ for all characters $\chi$ of $H$. Then
$\widehat{H}_{\psi}^{0}(G,M)=0$ for all characters
$\psi$ of $G$ and (consequently)
$$[M\e]=\prod_{\psi\in\widehat{G}}[M^{\e\psi}\e]=
\prod_{\psi\in\widehat{G}}[\e\varepsilon_{\psi} M\e].$$
\endproclaim
{\it Proof}. The stated formula will follow from Theorem 2.2 once we prove the first assertion. Let $\psi$ be any character of $G$ and let $M_{\e\overline{\psi}}$ be the twist of $M$ by $\overline{\psi}$ (see \S 1). As an $R[H]$-module, $M_{\e\overline{\psi}}$ is naturally isomorphic to $M_{\e\overline{\chi}}$, where $\overline{\chi}=\overline{\psi}|_{H}$ is the restriction of $\overline{\psi}$ to $H$. Therefore
$\widehat{H}^{0}(H,M_{\e\overline{\psi}})\simeq \widehat{H}^{0}(H,M_{\e\overline{\chi}})\simeq \widehat{H}^{0}_{\chi}(H,M)=0$. Now Theorem 4.4 shows that  $\widehat{H}^{0}(G,M_{\e\overline{\psi}})$ is zero, so $\widehat{H}^{0}_{\psi}(G,M)$
is zero and the result follows.\qed 
\smallpagebreak

{\it Remarks}. (a) Note that the number of conditions to be checked in order to apply Theorem 4.10 is {\it independent of $n$}.
See the remarks following the proof of Theorem 4.6.

(b) By Corollary 4.8 and the proof of the theorem, the theorem applies in fact to those finite $R[G\e]$-modules $M$ for which
the twists $M_{\xi^{i}}$, where $\xi$ is
a fixed generator of $\widehat{G}$ and $i=0, 1, \dots p-1$, are {\it cohomologically trivial}.
\smallpagebreak

When $p=2$, Theorem 4.10 leads to the following satisfying statement.
\proclaim{Corollary 4.11} Let $G$ be a cyclic group of order $2^n$, $n\geq 1$, and let $\tau$ denote
the unique element of $G$ of order 2. Let $M$ be a
finite $R[G\e]$-module such that $M^{+}=\{m\in M\:
\tau\e m=m\}=(1+\tau\e)M$.
Then
$$[M\e]=\prod_{\psi\in\widehat{G}}[M^{\e\psi}\e]
=\prod_{\psi\in\widehat{G}}[\e\varepsilon_{\psi}\be M\e].$$
\endproclaim
{\it Proof}. This will follow from Theorem 4.10 once
we check that $\widehat{H}_{\chi}^{0}(H,M)=0$, where $H=\langle\tau\rangle$ and $\chi$ is the
nontrivial character of $H$. But 
$$\widehat{H}_{\chi}^{0}(H,M)=
\sideset_{(1+\tau)}\and\to{\be M}/
(1-\tau)M=\widehat{H}^{-1}(H,M),$$
and the latter group has the same order as $\widehat{H}^{0}(H,M)=M^{+}/(1+\tau)\e M=0$.\qed

We now generalize Theorem 4.10 to arbitrary finite abelian groups by combining it with Theorem 3.1.
\proclaim{Theorem 4.12} Let $G$ be any finite abelian group and let $M$ be a finite $R[G\e]$-module, where $R$ is large enough for $G$. Suppose $G=K_{1}\times\dots
\times K_{r}$ is a decomposition of $G$ as a direct product of cyclic groups of $p_{i}$-power order, where 
$p_{1},p_{2},\dots,p_{r}$ are (not necessarily distinct) primes.
For each $0\leq i\leq r-1$ write $G_{i}=K_{1}\times\dots\times K_{i-1}$
($G_{0}=\{0\}$). Further,
for each $1\leq i\leq r$, let $H_{i}$ denote
the unique subgroup of $K_{i}$ of order $p_{i}$. Assume that $\widehat{H}_{\chi}^{0}(H_{i},M^{\e\psi})=0$ for all characters $\chi$ of $H_{i}$ and
$\psi$ of $G_{i-1}$, for $i=1,2,\dots, r$.
Then
$$[M\e]=\prod_{\psi\in\widehat{G}}[M^{\e\psi}\e]
=\prod_{\psi\in\widehat{G}}[\e\varepsilon_{\psi}\be M\e].$$
\endproclaim
{\it Proof}. This follows at once from Theorems 3.1 and 4.10.\qed
\smallpagebreak

{\it Remark}. The above theorem cannot be considered satisfactory because of its rather stringent hypotheses. We have stated it only for completeness and will not give applications of it.

\heading 5. Applications to abelian varieties\endheading

Let $F$ be a global field, i.e., $F$ is a finite extension of ${\Bbb Q}$ (the ``number field case") or is finitely generated and of
transcendence degree 1 over a finite field (the ``function field case"). Let $A$ be an abelian variety defined over $F$ and let $K/F$ be a finite Galois extension with Galois group $G$. We will write $A_{F}$ (resp. $A_{K}$) for the abelian variety $A$ regarded as an abelian variety over $F$ (resp. $K$), and $A^{\roman{t}}$ will denote the abelian variety dual to $A$.
Further, for any prime $w$ of $K$, we will
write $G_{w}$ for ${\roman{Gal}}(K_{w}/F_{v})$, 
where $v$ is the prime of $F$ lying below $w$, and will identify this group with the
decomposition group of $w$ in $G$.

\proclaim{Proposition 5.1} Let $w$ be any prime of $K$, let $v$ be the prime of $F$ lying below $w$ and let $L/F_{v}$ be a subextension of $K_{w}/F_{v}$ such that $H=
\roman{Gal}(K_{w}/L)$ is cyclic. Then
the Herbrand quotient of the $H$-module
$A(K_{w})$ is 1.
\endproclaim
{\it Proof}. Let $f\: A\ra A^{\roman{t}}$ be any isogeny and write $A_{f}$ for the kernel of $f$. Then there is an exact sequence
$$0\ra A_{f}(\overline{K}_{w})
\ra A(\overline{K}_{w})\ra A^{\roman{t}}(\overline{K}_{w})\ra 0, $$
where $\overline{K}_{w}$ is a separable algebraic closure of $K_{w}$. Taking $\roman{Gal}(\overline{K}_{w}/K_{w})$-invariants of the above exact sequence, we
conclude that there exists an
$H$-module homomorphism $A(K_{w})\ra A^{\roman{t}}(K_{w})$ having a finite kernel and cokernel. Consequently $h(A(K_{w}))=h(A^{\roman{t}}(K_{w}))$. On the other hand local duality [9, I.3.4, 3.7; III.7.8] (see also [8, 4.2]) shows that $h(A(K_{w}))\e\e h(A^{\roman{t}}(K_{w}))=1$. Thus $h(A(K_{w}))^2=1$, whence $h(A(K_{w}))=1$.\qed
\smallpagebreak

{\it Remark}. There is an alternative proof
of Proposition 5.1 in the case where $K_{w}$ is non-archimedean of characteristic zero. This proof uses a well-known theorem of Mattuck and therefore depends on the theory of the logarithm. See [12, \S 4, (14)].
\smallpagebreak

\proclaim{Corollary 5.2} Let $v$ be a real prime of $F$.  Then $A(F_{v})$ is connected
if and only $A^{\roman{t}}(F_{v})$ is connected.
\endproclaim
{\it Proof}. By [9, I.3.7], $A(F_{v})$ is connected if and only if $\widehat{H}^{0}(F_{v},A)=0$ or, equivalently, if and only
if $H^{1}(F_{v},A^{\roman{t}})=0$ (by duality). The proposition (applied to some totally imaginary finite Galois extension of $F$) shows that the latter is equivalent to the vanishing of $\widehat{H}^{0}(F_{v},A^{\roman{t}})$, i.e., to the connectedness of $A^{\roman{t}}(F_{v})$.\qed
\smallpagebreak

{\it Remark.} Corollary 5.2 was already known to Yu.Zarhin in 1972. See [14, \S 3].
\smallpagebreak

\proclaim{Theorem 5.3} Let $v$ be any prime of $F$ and let $w$ be a fixed prime of $K$ lying above $v$. Assume that $A(L)=N_{K_{w}/L}A(K_{w})$ for all subextensions $L/F_{v}$
of $K_{w}/F_{v}$ of prime index, i.e., such that $[K_{w}:L]$ is prime. Then the $G_{w}$-module $A(K_{w})$ is cohomologically trivial.
\endproclaim
{\it Proof}. The theorem follows from Proposition 5.1 and Corollary 4.8.
\qed
\smallpagebreak

{\it Remark}. With the notations of the theorem,
suppose that $K_{w}/L$ is ramified of prime
degree $p$. Write $l$ for the residue field of $L$. Further, let $\Cal A$ denote the N\'eron model of $A_{L}$ and $\widehat{{\Cal A}}$ the formal completion of $\Cal A$ along its zero section. Then
there is an exact sequence
$$A(l)_{p}\ra\widehat{H}^{0}(H,\widehat{{\Cal A}}(K_{w}))\ra
\widehat{H}^{0}(H,A(K_{w}))\ra A(l)/p A(l)\ra 0,$$
in which $A(l)_{p}$ denotes the $p$-torsion subgroup of $A(l)$ (see [8, Corollary 4.6]). It follows that the vanishing of $\widehat{H}^{0}(H,A(K_{w}))$ is equivalent to the vanishing of both $A(l)_{p}$ and
$\widehat{H}^{0}(H,\widehat{{\Cal A}}(K_{w}))$. The vanishing of $A(l)_{p}$ may occur with some frequency,
while it may be possible to find simple conditions under
which $\widehat{H}^{0}(H,\widehat{{\Cal A}}(K_{w}))$ is zero, perhaps using appropriate variants of the methods of [8, \S 4].
However, we have not looked into this problem closely and are therefore unable to comment further on it at the present time.
\smallpagebreak

Next we recall the main theorem of [6]. We will
write $T$ for the set of primes of $F$ which ramify in $K/F$ or where $A$ has bad reduction. 
\proclaim{Theorem 5.4} Notations being as above, suppose that the following conditions hold.
\roster
\item"(i)" The Tate-Shafarevich group ${\cyr W}(A_{K})$ of $A_{K}$ is finite.
\item"(ii)" $\widehat{H}^{i}(G,A(K))=
\widehat{H}^{i}(G,A^{\roman{t}}(K))=
0$ for all integers $i$.
\item"(iii)" $A(F_{v})$ is connected for all real primes $v$ of $F$.
\endroster
Then 
$$\left[{\cyr W}(A_{K})^{G}\e\right]=[{\cyr W}(A_{F})]\cdot\prod_{v\in T}[H^{1}(G_{w},A(K_{w}))]$$
and
$$[H^{1}(G,{\cyr W}(A_{K}))]=\prod_{v\in T}[H^{2}(G_{w},A(K_{w}))],$$
where, for each prime $v\in T$, $w$ is
a fixed prime of $K$ lying above $v$.
\endproclaim
{\it Proof}. See [6], Theorem 4.4.\qed
\smallpagebreak

{\it Remarks}. (a) The above result, which was established in [6] for number fields only, is in fact valid for any global field. This is so because [9], which was the main reference for [6], covers both the number field and function field cases (one only needs to supplement some of the references made in [6] to results from Chapter I of [9] with
references to Chapter III of the same book).

(b) As pointed out by the referee of [6], the conditions of the theorem ``[seem] rather stringent but hold in fact quite often". Aljadeff's theorem ultimately
tells us why this is so. Consider for example the case of a $p$-group $G$, where $p$ is a prime number. Then condition (ii) of the theorem is equivalent to the cohomological triviality of both
$A(K)$ and $A^{\roman{t}}(K)$. By Theorem 4.6, the
latter is equivalent to the vanishing of
$\widehat{H}^{i}(H,A(K))$ and
$\widehat{H}^{i}(H,A^{\roman{t}}(K))$ for $i=-1, 0$ and for all subgroups $H$ of $G$ of order $p$\footnote"*"{Furthermore, observe that if there exists an isogeny $f\: A\ra A^{\roman{t}}$ of degree prime to $p$, then
$\widehat{H}^{i}(H,A^{\roman{t}}(K))=0$ for all $i$ if and only if $\widehat{H}^{i}(H,A(K))=0$ for all $i$. See the proof of Proposition 5.1 above.}. If, furthermore, $G$ is cyclic (so that $G$ has a unique subgroup of order $p$), then the above conditions do not seem stringent at all.

(c) Condition (iii) is vacuous if $F$ has no real primes. Furthermore, Corollary 5.2 shows that it is equivalent to condition (B) of [6].
\smallpagebreak

We assume now, for simplicity, that $K/F$ is a {\it cyclic} extension\footnote{It is possible to give versions of the results of this section assuming only that the extension $K/F$ is abelian, but the statements are too complicated and we omit them.}.

\proclaim{Corollary 5.5} Assume that $G$ is cyclic. Then, under the assumptions of Theorem 5.4,
$$\big[{\cyr W}(A_{K})^{G}\big]=[{\cyr W}(A_{F})]\e\e[\widehat{H}^{0}(G,{\cyr W}(A_{K}))]$$
and
$$[\widehat{H}^{0}(G,{\cyr W}(A_{K}))]=
\prod_{v\in T}[\widehat{H}^{0}(G_{w},A(K_{w}))].$$
\endproclaim
{\it Proof}. This follows at once from Theorem 5.4, using Proposition 5.1, the periodicity of the cohomology of cyclic groups, and the fact that the Herbrand quotient of a finite module is 1.\qed

We now write $n$ for the order of $G$ and assume that $A$ has {\it complex multiplication} by the ring of integers $R={\Bbb Z}[\zeta_{n}]$ of ${\Bbb Q}(\zeta_{n})$. Then ${\cyr W}(A_{K})$ is an $R[G\e]$-module in a natural way, and we may therefore apply to it the results of the preceding sections. For each character
$\chi$ of $G$, we will write $A^{\chi}$ for the $\chi$-twist of $A$ \footnote{
This is a standard notation for the $\chi$-twist of an abelian variety. We have adopted it in spite of the fact that some readers may be confused by it, in view of the notations introduced earlier. It may help clarify matters to note that $A^{\chi}(K)$ is an $R[G\e]$-module which is isomorphic to the twisted $R[G\e]$-module $A(K)_{\chi}$ defined in \S 1, whereas $A(K)^{\chi}$ (which we primarily regard as an $R$-module) is isomorphic to the $R$-module $A^{\overline{\chi}}(F)$, where $\overline{\chi}$ is, as usual, the inverse of $\chi$.}(see [10, \S 2]). 
Then there are isomorphisms
$${\cyr W}(A_{K})^{\chi}\simeq({\cyr W}(A_{K})_{\e\overline{\chi}})^{G}
\simeq{\cyr W}(A^{\overline{\chi}}_{\be\be K})^{G}.$$
The next corollary results from applying
Corollary 5.5 to the twisted abelian variety $A^{\overline{\chi}}$.
\proclaim{Corollary 5.6} Let $G$ be cyclic and let $\chi$ be a character of $G$. Assume that the conditions of Theorem 5.4 hold for the twisted 
abelian variety $A_{\be\be K}^{\chi}$. Then
$$[{\cyr W}(A_{K})^{\chi}\e]=[{\cyr W}(A^{\overline{\chi}}_{\be\be F})]\e\e[\widehat{H}_{\chi}
^{0}(G,{\cyr W}(A_{K}))]$$
and
$$[\widehat{H}_{\chi}^{0}(G,{\cyr W}(A_{K}))]=
\prod_{v\in T_{\e\overline{\chi}}}[\widehat{H}_{
\chi_{_{w}}}^{0}\be\be(G_{w},A(K_{w}))],$$
where $T_{\e\overline{\chi}}$ is the set of primes of $F$ that ramify in $K/F$ or where $A^{\overline{\chi}}_{\be\be F}$ has bad reduction and $\chi_{_{w}}$ is the restriction of $\chi$ to $G_{w}$.
\endproclaim

We now combine the preceding corollary with Theorem 2.2 to obtain the main result of this section.

\proclaim{Theorem 5.7} Let $K/F$ be a cyclic Galois
extension of global fields with Galois group $G$ of order $n$ and let $A$ be an abelian variety over
$F$ with complex multiplication by the ring of integers of ${\Bbb Q}(\zeta_{n})$. Assume that
the following conditions hold.
\roster
\item"(i)" The Tate-Shafarevich group ${\cyr W}(A_{K})$ of $A_{K}$ is finite.
\item"(ii)" $\widehat{H}_{\chi}^{i}(G,A(K))=
\widehat{H}_{\chi}^{i}(G,A^{\roman{t}}(K))=
0$ for all $i$ and all $\chi\in
\widehat{G}$.
\item"(iii)" $A(F_{v})$ is connected for all real primes $v$ of $F$.
\endroster
Then
$$[{\cyr W}(A_{K})]=\prod_{\chi\in\widehat{G}}
[{\cyr W}(A^{\chi}_{\be\be F})]\cdot\prod_{\chi\in\widehat{G}}[S_{\chi}(
{\cyr W}(A_{K}))],$$
where, for each $\chi\in
\widehat{G}$, $S_{\chi}({\cyr W}(A_{K}))$ is the submodule of $\widehat{H}_{\chi}^{0}(G,{\cyr W}(A_{K}))$ defined in Section 2, formula (6). Further,
the order of $S_{\chi}({\cyr W}(A_{K}))$ divides the product
$$\prod_{v\in T_{\e\overline{\chi}}}[\widehat{H}_{\chi_{_{w}}}^{0}\be\be(G_{w},A(K_{w}))],$$
where $T_{\e\overline{\chi}}$ and $\chi_{_{w}}$ are as in the statement of Corollary 5.6.
\endproclaim
\smallpagebreak

{\it Remarks}. (a) Since $A$ and $A^{\chi}$ are isomorphic over $K$, the finiteness of ${\cyr W}(A^{\chi}_{K})$ is equivalent to that of ${\cyr W}(A_{K})$. Thus condition (i) of Theorem 5.7 implies that condition (i) of Theorem 5.4 holds for all twists of $A$.

(b) Similarly, condition (iii) of Theorem 5.7 suffices to ensure that condition (iii) of Theorem 5.4 holds for all twists of $A$. Indeed, let $v$ be a real prime of $F$ and let $K$ be a totally imaginary extension of $F$. By [9, I.3.7], $A^{\chi}(F_{v})$ is connected if and only if $\widehat{H}^{0}(G_{w},A^{\chi}(K_{w}))\simeq\widehat{H}_{\e\overline{\chi}_{_w}}^{0}\be\be(G_{w},A(K_{w}))$ is zero. On the other hand, since $K_{w}/F_{v}$ is a quadratic extension,
$$[\widehat{H}_{\e\overline{\chi}_{_{w}}}^{0}\be\be
(G_{w},A(K_{w}))]=[\widehat{H}^{-1}(G_{w},A(K_{w}))]=
[\widehat{H}^{0}(G_{w},A(K_{w}))],$$
by Proposition 5.1 and the proof of Corollary 4.11. Thus $A^{\chi}(F_{v})$ is connected if and only if $\widehat{H}^{0}(G_{w},A(K_{w}))$ is zero, i.e., if and only if $A(F_{v})$ is connected.

(c) If $G$ is a cyclic $p$-group, where $p$ is a prime, and $H$ is the unique subgroup
of $G$ of order $p$, then condition (ii) of Theorem 5.7 is equivalent to:
\roster
\item"(ii)$^{\prime}$" $\widehat{H}_{\chi}^{i}(H,A(K))=\widehat{H}_{\chi}
^{i}(H,A^{\roman{t}}(K))=0$ for $i=-1, 0$ and all $\chi\in\widehat{H}$.
\endroster
See the remark following the statement of Theorem 5.4.

(d) Theorem 5.7 generalizes Corollary 4.6 of [6]. As explained in the introduction, the search for such a generalization led to the writing of this paper.

\smallpagebreak

\heading 6. Applications to class groups of number fields\endheading

If $G$ is a finite group and $M$ is a finite
$\Bbb Z[G\e]$-module, then the results of \S 4 do not immediately apply to $M$ because $\Bbb Z$ is not large enough for $G$. Therefore we need to extend scalars, i.e., consider the ring $R
=\Bbb Z[\zeta_{e}]$, where $e$ is the exponent of $G$ and $\zeta_{e}$ is a fixed (complex) $e$-th root of unity. Define
$$\overline{M}= M\otimes_{\Bbb Z}R.$$
Clearly, $\overline{M}$ is an $R[G\e]$-module of order $[M\e]^{\varphi(e)}$, where $\varphi$ is Euler's function. More precisely, let $\Cal B=\{\zeta_{e}^{i}\: 0\leq i\leq \varphi(e)-1\}$. Then the elements of $\overline{M}$ may be regarded as formal linear combinations of elements of $M$ with coefficients in $\Cal B$, i.e., any $x\in\overline{M}$ can be written in a unique way in the form
$$x=\sum_{i=0}^{\varphi(e)-1}\be\be m_{i}\e\zeta_{e}^{i}\e ,\tag 7$$

We now apply Theorem 3.1 to the $R[G\e]$-module $\overline{C}_{K}=C_{K}\otimes_{\Bbb Z}R$ and obtain the following result.
\proclaim{Theorem 6.1} Let $K/F$ be an abelian extension of number fields with Galois group $G$ of exponent $e$. Let $G=H_1\times\dots\times H_r$ ($r\geq 1$) be a decomposition of $G$ as a direct product of cyclic groups. For $0\leq i\leq r-1$, let $G_{i}=H_1\times\dots\times H_{i}$, where $G_{0}$ is defined to be $\{0\}$. Then
$$(h_{K}/h_{F})^{\varphi(e)}\cdot\prod_{i=1}^{r}\prod_
{\Sb \chi_{_i}\be\be\in\widehat{H}_{i}\\
\psi\in\widehat{G}_{i-1}\endSb}
\ng\ng\left[\widehat{H}^0_{\chi_{_i}}(H_i,\overline{C}_{K}^{\e\psi}\e)/S_{\chi
_{_i}}(\overline{C}_{K}^{\e\psi}\e)\right]=\prod_{
\Sb \chi\in\widehat{G}\\\chi\neq\chi^{0}\endSb}
\e[\overline{C}_{K}^{\e\chi}\e].$$
\endproclaim
\smallpagebreak

Admittedly the above formula is rather complicated, but is not without interesting applications, as we now show. One immediate observation is that the product which appears on the left-hand side of the above formula is an integer which is divisible {\it only by the primes that divide $e$}. This comes from the fact that multiplication by $e$
annihilates each of the groups
$\widehat{H}^0_{\chi_{_i}}(H_i,\overline{C}_{K}^
{\e\psi}\e)$. Further, writing each of the factors on the right-hand side in the form $[\varepsilon_{\chi}\overline{C}_{K}\e]
[\widehat{H}^0_{\chi}(G,\overline{C}_{K})]$, we see at once that there exists a positive rational number $r$, whose numerator and denominator are divisible only by primes that divide $e$, such that
$$r\cdot(h_{K}/h_{F})^{\varphi(e)}=\prod_{
\chi\neq\chi^{0}}
[\varepsilon_{\chi}\overline{C}_{K}\e].\tag 8$$
Thus, for example, if $G$ is a $p$-elementary abelian group, where $p$ is a prime, then the formula of the theorem is an identity of the type
$$p^{t}\cdot(h_{K}/h_{F})^{\e p-1}
=\prod_{\chi\neq\chi^{0}}
\e[\varepsilon_{\chi}\overline{C}_{K}],$$
where $t$ is an integer (which may be positive, negative or zero). On the probable value of $t$, see below. When $p=2$ the situation is particularly simple, because in this case there is no need to extend scalars ($R=\Bbb Z$). We have

\proclaim{Corollary 6.2} Let $K/F$ be a finite abelian
extension of exponent 2. Then there exists an integer $t$ such that
$$2^{t}\cdot (h_{K}/h_{F})=\prod_{\chi\neq\chi^{0}}
\e[\varepsilon_{\chi}C_{K}],$$
where the product extends over all non-trivial characters of $G$.
\endproclaim

Regarding the formula of the corollary,
the factors $[\varepsilon_{\chi}C_{K}]$ which appear on the right-hand side are (or seem to be) related to the class numbers of the various subextensions of $K/F$. For example.
\smallpagebreak

{\it Example}. ( See [13, Theorem 10.10, p.191])
Let $K=\Bbb Q(\sqrt{d_{1}},\sqrt{d_{2}})$ be a biquadratic extension of the rational field, where $d_{1}$ and $d_{2}$ are squarefree integers. Let
$\chi_{j}$ ($j=1,2,3$) denote the nontrivial
characters of $G=\roman{Gal}(K/F)$, and assume that the numbering has been chosen so that
$L_{j}=\roman{Fix}(\krn\chi_{j})=\Bbb Q(\sqrt{d_{j}})$, where $d_{3}=d_{1}d_{2}$. Write
$$\align\{1,\tau\}&=\roman{Gal}(K/L_{1})\\
\{1,\sigma\}&=\roman{Gal}(K/L_{2})\\
\{1,\sigma\tau\}&=\roman{Gal}(K/L_{3}),\endalign$$
and set $\varepsilon_{\chi_{j}}=\varepsilon_{j}$
($j=1,2,3$). Then $\varepsilon_{1}=(1-\sigma)(1+\tau)$, $\varepsilon_{2}=(1-\tau)(1+\sigma)$ and $\varepsilon_{3}=(1-\tau)(1-\sigma)$. We have
$\varepsilon_{1}C_{K}=
(1-\sigma)N_{K/L_{1}}C_{K}$ and, moreover, $\sigma$ acts on $N_{K/L_{1}}C_{K}$ as multiplication by $-1$ since
$$(1+\sigma)N_{K/L_{1}}C_{K}=N_{L_{1}/\Bbb Q}
N_{K/L_{1}}C_{K}=0.$$
It follows that $[\varepsilon_{1}C_{K}]$ differs from
$h_{1}\overset{\text{def}}\to{=}[C_{L_{1}}]$ by a power of 2. Similarly,
$[\varepsilon_{2}C_{K}]$ differs from $h_{2}=[C_{L_{2}}]$ by a power of 2. On the other hand $\varepsilon_{3}C_{K}=
(1-\tau)(1+\sigma\tau)C_{K}=(1-\tau)N_{K/L_{3}}C_{K}$,
and we conclude as before that
$[\varepsilon_{3}C_{K}]$ differs from $h_{3}=[C_{L_{3}}]$ by a power of 2. Summarizing, there exists an integer $t$ such that
$$2^{t}\cdot h_{K}=h_{1}\e h_{2}\e h_{3}.$$
\medpagebreak

{\it Remarks.} (a) The formula of the example is certainly not new. A more precise version of it [5, Theorem 74, p.320] follows from the classical Brauer relations [Ibid., Theorem 73, p.315]. Incidentally, it looks like an interesting problem to compare the class number formulas of this section with those obtained by R.Brauer through analytical means. Since the latter involve
the regulators of a field and its various subfields, while the formulas of this section do not, it seems
likely that some non-trivial relations among
regulators will emerge, generalizing [Op.Cit., formula (7.27), p.320]. Regarding this last remark, some readers may be surprised by the fact that the units of the field do not seem to play a role in our formulas. In fact, the units lurk in the background. See the next remark.

(b) To determine the exact value of $r$ in (8), at least under the assumption that $K/F$ is cyclic, one would proceed as follows. It is known that the ideal class group of a number field $K$ is the Tate-Shafarevich group of the units of $K$, i.e.,
$$C_{K}\simeq\krn\big[H^{1}(K, U_{K})\ra\prod H^{1}(K_{w},
U_{w})\big].$$
See [4]. Thus, by using the methods of [6], it should not be difficult to obtain a formula for the order of $\widehat{H}^{0}(G, C_{K})$
similar to the identity $[\widehat{H}^{0}(G,{\cyr W}(A_{K}))]=\prod_{v\in T}[\widehat{H}^{0}(G_{w},A(K_{w}))]$ of Corollary 5.5,
at least under the additional assumption that the analog of  condition (ii) of Theorem 5.4, namely $\widehat{H}^{0}(G,U_{K})=0$, holds. Without this assumption, a formula for $[\widehat{H}^{0}(G, C_{K})]$ should involve the ramification indices
$e(w/v)=[\widehat{H}^{0}(G_{w},U_{w})]$ (where $v$ is the prime of $F$ lying below $w$), the order of $\widehat{H}^{0}(G,U_{K})$, and possibly some other familiar constants. For a clear indication that this is so see [11], where the author essentially computes a ``toric analogue" of  $[\widehat{H}^{0}(G, C_{K})]$. The above comments, although not immediately applicable to the problem of computing 
$[\widehat{H}^0_{\chi_{_i}}(H_i,\overline{C}_{K}^
{\e\psi}\e)]$ for non-trivial characters $\chi_{i}$
and $\psi$, should however give an indication as to what to expect regarding the value of $r$ in (8).
\smallpagebreak

The next theorem generalizes a well-known consequence of class field theory.

\proclaim{Theorem 6.3} Let $K/F$ be a finite Galois extension of number fields. Suppose that for every subextension of prime index $L/F$ of $K/F$ the extension $K/L$ ramifies at some prime. Then the norm map $N_{K/L}\: C_{K}\ra C_{L}$ is surjective for each subextension $L/F$ of $K/F$. In particular, $h_{L}$ divides $h_{K}$ for every subextension $L/F$ of $K/F$.
\endproclaim
{\it Proof}. The hypothesis implies that $K\cap H_{L}=L$ for every $L/F\subset K/F$ of prime index, where $H_{L}$ is the Hilbert class field of $L$. A well-known consequence of class field theory [7, Lemma, p. 83] then shows that
$N_{K/L}\: C_{K}\ra C_{L}$ is surjective for all subextensions $L/F$ of $K/F$ of prime index. Theorem 4.4 now implies that $N_{K/L}$ is surjective for all subextensions $L/F$ of $K/F$. The last assertion of the theorem is clear.\qed
 
\proclaim{Theorem 6.4} Let $K/F$ be a cyclic Galois extension of number fields with Galois group $G$ of order $n$. Assume that $K/F$ satisfies the conditions of
Theorem 6.3 above. Then
$$(h_{K}/h_{F})^{\varphi(n)}=\prod_{\Sb \chi\in\widehat{G}\\\chi\neq\chi^{0}\endSb}
\e[S_{\chi}(\overline{C}_{K})]\cdot
\prod_{\Sb \chi\in\widehat{G}\\\chi\neq\chi^{0}\endSb}
\e[\varepsilon_{\chi}\overline{C}_{K}],$$
where the integers $[S_{\chi}(\overline{C}_{K})]$ (which were defined in Section 2, formula (6)) are divisible only by the primes that divide $n$.
\endproclaim
{\it Proof.} This follows from Theorem 2.2.\qed
\smallpagebreak

{\it Remarks}. (a) As noted earlier, the factors
$[\varepsilon_{\chi}\overline{C}_{K}]$ which appear in the formula of the theorem seem to be related to the class numbers of the various subextensions of $K/F$. In particular if $n=p$ is a prime the identity of the theorem should reduce to a trivial factorization of the type $(h_{K}/h_{F})^{p-1}=
(h_{K}/h_{F})\cdots (h_{K}/h_{F})$ ($p-\be 1$ factors). In general, the computation of $[\varepsilon_{\chi}\overline{C}_{K}]$ (or $[\overline{C}_{K}^{\e\e\chi}]$)
in terms of class numbers of subextensions of $K/F$ seems to be a problem in linear algebra. For example, if $\tau$ is a fixed generator of $G$ and $\chi$ is given by $\chi(\tau)=\zeta_{n}$, then there is a natural isomorphism
$$\overline{C}_{\be K}^{\e\e\chi}\simeq\{\bold{m}\in C_{K}^{\e \varphi(n)}\:\tau\bold{m}=A\e\bold{m}\},$$
where $A$ is the {\it companion matrix}\,\footnote{Or the {\it transpose} of the companion matrix, depending on which definition of ``companion matrix" one adopts.} of the $n$-th cyclotomic polynomial $\Phi_{n}(x)$ (this follows from  (7). Incidentally, if $n$ is a prime then the above ``eigenspace" is naturally isomorphic to $C_{K}$). However, we do not yet know if $\overline{C}_{\be K}^{\e\e\chi^{i}}$ ($2\leq i\leq n-1$) admits a similar description. Whatever the case may be, it appears at this time that the identity
$$N_{K/L}=\prod_{\Sb d\mid[K:L]\\d>1\endSb}
\Phi_{d}\big(\tau^{[L:F]}\big)$$
will play a role in the computation of $[\overline{C}_{\be K}^{\e\e\chi^{i}}]$ in terms of the class numbers of the various subextensions of $K/F$.

(b) Let $p$ be an odd prime and let
$K=\Bbb Q(\zeta_{p})^{+}$ be the maximal real subfield
of $\Bbb Q(\zeta_{p})$. We write
$h^{+}$ for the class number of $K$. For $F=\Bbb Q$, the formula of the theorem becomes
$$(h^{+})^{\varphi\left(\frac{p-1}{2}\right)}=r\cdot
\prod_{\chi\neq\chi^{0}}\e[\varepsilon_{\chi}\overline{C}_{K}],$$
where $r$ is an integer which is divisible only by the primes that divide $(p-1)/2$. In particular,
$p$ divides $h^{+}$ if and only if $p$ divides $[\varepsilon_{\chi}\overline{C}_{K}]$ for some
(non-trivial) character $\chi$ of $G$, i.e., Vandiver's conjecture holds for $p$ if and only if $p$ does not divide $[\varepsilon_{\chi}\overline{C}_{K}]$
for each character $\chi$ of $G$. Concerning this conjecture, the results of this paper (cf. previous remark) seem to indicate that the following statement is true: $p$ does not divide $h^{+}$ if and only if $p$ does not divide $h_{L}$ for every subextension $L/\Bbb Q$ of $K/\Bbb Q$ of {\it prime} degree. More generally, the results of this section seem to suggest that for a large class of abelian extensions of a given number field (perhaps that which consists of the fields satisfying the conditions of Theorem 6.3), the class numbers of its subextensions of prime degree to a large extent determine
the class numbers of {\it all} its subextensions.
\medpagebreak

We conclude this paper with the following
statement.
\proclaim{Theorem 6.5} Let $K/F$ be a cyclic Galois extension with Galois group $G$ of order $2^{n}$, where $n\geq 1$. Assume that $K/F$ ramifies at some prime. Then
$$(h_{K}/h_{F})^{2^{n-1}}=\prod_
{\chi\neq\chi^{0}}[\e\varepsilon_{\chi}\be
\overline{C}_{K}].$$
\endproclaim
{\it Proof.} The hypothesis implies that $K/L$ ramifies at some prime, where $L/F$ is the unique subextension of $K/F$ of index 2. Therefore, by Theorem 6.4, $N_{K/L}$ is surjective. The theorem now
follows from Corollary 4.11.\qed

\enddocument